%
%
%

\def\LaTeX{%
  \let\Begin\begin
  \let\End\end
  \let\salta\relax
  \let\finqui\relax
  \let\futuro\relax}

\def\UK{\def\our{our}\let\sz s}
\def\USA{\def\our{or}\let\sz z}



\LaTeX

\USA


\salta

\documentclass[twoside,12pt]{article}
\setlength{\textheight}{24cm}
\setlength{\textwidth}{16cm}
\setlength{\oddsidemargin}{2mm}
\setlength{\evensidemargin}{2mm}
\setlength{\topmargin}{-15mm}
\parskip2mm


\usepackage{amsmath}
\usepackage{amsthm}
\usepackage{amssymb}
\usepackage[mathcal]{euscript}

\usepackage[usenames,dvipsnames]{color}
%
%
%
\def\gianni{\marginpar{{$\leftarrow$}}\color{green}}
\def\piecol{\color{red}}
\def\pier{\color{red}}
\def\ppg{\color{blue}}           
\def\giannipier{\color{magenta}}     
%
%
\let\gianni\relax
\let\piecol\relax
\let\pier\relax
\let\ppg\relax
\let\giannipier\relax


\bibliographystyle{plain}


%

\finqui

\def\Beq{\Begin{equation}}
\def\Eeq{\End{equation}}
\def\Bsist{\Begin{eqnarray}}
\def\Esist{\End{eqnarray}}

\def\Bthm{\Begin{theorem}}
\def\Ethm{\End{theorem}}
\def\Blem{\Begin{lemma}}
\def\Elem{\End{lemma}}

\def\Brem{\Begin{remark}\rm}
\def\Erem{\End{remark}}

\def\Bdim{\Begin{proof}}
\def\Edim{\End{proof}}
\let\non\nonumber




\def\step #1 \par{\medskip\noindent{\bf #1.}\quad}


\def\holder{H\"older}
\def\aand{\quad\hbox{and}\quad}

\def\lhs{left-hand side}
\def\rhs{right-hand side}

\def\wk{well-known}


\def\organiz{organi\sz}

\def\regulariz{regulari\sz}
\def\summariz{summari\sz}


\def\multibold #1{\def\arg{#1}%
  \ifx\arg\pto \let\next\relax
  \else
  \def\next{\expandafter
    \def\csname #1#1#1\endcsname{{\bf #1}}%
    \multibold}%
  \fi \next}

\def\pto{.}

\def\multical #1{\def\arg{#1}%
  \ifx\arg\pto \let\next\relax
  \else
  \def\next{\expandafter
    \def\csname cal#1\endcsname{{\cal #1}}%
    \multical}%
  \fi \next}


\def\multimathop #1 {\def\arg{#1}%
  \ifx\arg\pto \let\next\relax
  \else
  \def\next{\expandafter
    \def\csname #1\endcsname{\mathop{\rm #1}\nolimits}%
    \multimathop}%
  \fi \next}

\multibold
qwertyuiopasdfghjklzxcvbnmQWERTYUIOPASDFGHJKLZXCVBNM.

\multical
QWERTYUIOPASDFGHJKLZXCVBNM.

\multimathop
dist div dom meas sign supp .


\def\accorpa #1#2{\eqref{#1}--\eqref{#2}}
\def\Accorpa #1#2 #3 {\gdef #1{\eqref{#2}--\eqref{#3}}%
  \wlog{}\wlog{\string #1 -> #2 - #3}\wlog{}}


\def\tonde #1{\left(#1\right)}

\def\graffe #1{\mathopen\{#1\mathclose\}}

\def\<#1>{\mathopen\langle #1\mathclose\rangle}
\def\norma #1{\mathopen \| #1\mathclose \|}

\def\iot {\int_0^t}

\def\iO{\int_\Omega}
\def\intQt{\int_{Q_t}}

\def\dt{\partial_t}
\def\dn{\partial_\nu}
\def\ds{\,ds}

\def\cpto{\,\cdot\,}

\def\checkmmode #1{\relax\ifmmode\hbox{#1}\else{#1}\fi}
\def\aeO{\checkmmode{a.e.\ in~$\Omega$}}
\def\aeQ{\checkmmode{a.e.\ in~$Q$}}
\def\aeS{\checkmmode{a.e.\ on~$\Sigma$}}

\def\aat{\checkmmode{for a.a.~$t\in(0,T)$}}


\def\erre{{\mathbb{R}}}




\def\genspazio #1#2#3#4#5{#1^{#2}(#5,#4;#3)}
\def\spazio #1#2#3{\genspazio {#1}{#2}{#3}T0}

\def\L {\spazio L}
\def\H {\spazio H}
\def\W {\spazio W}

\def\C #1#2{C^{#1}([0,T];#2)}


\def\Lx #1{L^{#1}(\Omega)}
\def\Hx #1{H^{#1}(\Omega)}
\def\Wx #1{W^{#1}(\Omega)}
\def\Cx #1{C^{#1}(\overline\Omega)}

\def\Luno{\Lx 1}
\def\Ldue{\Lx 2}
\def\Linfty{\Lx\infty}
\def\Lq{\Lx4}
\def\Huno{\Hx 1}
\def\Hdue{\Hx 2}


\def\LQ #1{L^{#1}(Q)}
\def\CQ {C^0(\overline Q)}


\let\theta\vartheta
\let\eps\varepsilon

\let\TeXchi\chi                         
\newbox\chibox
\setbox0 \hbox{\mathsurround0pt $\TeXchi$}
\setbox\chibox \hbox{\raise\dp0 \box 0 }
\def\chi{\copy\chibox}


\def\muz{\mu_0}
\def\rhoz{\rho_0}
\def\uz{u_0}

\def\normaV #1{\norma{#1}_V}
\def\normaH #1{\norma{#1}_H}
\def\normaW #1{\norma{#1}_W}

\def\kamurho{\kappa(\mu,\rho)}
\def\kapmurho{\kappa'(\mu,\rho)}
\def\kmin{\kappa_*}
\def\kmax{\kappa^*}

\def\coeff{1+2g(\rho)}
\def\coeffs{1+2g(\rho(s))}
\def\coefft{1+2g(\rho(t))}

\def\T{\calT_\tau}

\def\mut{\mu_\tau}
\def\rhot{\rho_\tau}
\def\xit{\xi_\tau}
\def\kappat{\bar\kappa}
\def\Ku{K_1}
\def\Kd{K_2}

\def\kamurhos{\kappa(\mu(s),\rho(s))}
\def\Kmurho{K(\mu,\rho)}
\def\Kumurho{\Ku(\mu,\rho)}
\def\Kdmurho{\Kd(\mu,\rho)}
\def\Kmurhot{K(\mu(t),\rho(t))}
\def\Kumurhot{\Ku(\mu(t),\rho(t))}
\def\Kmurhos{K(\mu(s),\rho(s))}

\def\coefftau{1+2g(\rhot)}

\def\coeffn{1+2g(\rhon)}
\def\coeffnt{1+2g(\rhon(t))}

\def\mun{\mu_n}
\def\rhon{\rho_n}
\def\xin{\xi_n}
\def\munmu{\mu_{n-1}}

\def\coeffn{1+2g(\rhon)}

\def\mumkp{(\mu-k)^+}

\def\iotmt{\int_0^{t-\tau}}
\def\intQtmt{\iotmt\!\!\!\iO}



\Begin{document}


\title{\bf Global existence for a strongly coupled 
Cahn-Hilliard system with viscosity%
\footnote{{\bf Acknowledgments.}\quad\rm P. Colli and G. Gilardi 
gratefully acknowledge the financial support of
the MIUR-PRIN Grant 2008ZKHAHN {\sl ``Phase transitions, 
hysteresis and multiscaling''} and of the IMATI of CNR in Pavia. 
The work of J. Sprekels was supported by the DFG Research Center 
{\sc Matheon} in Berlin.}}
\author{}
\date{}
\maketitle
\begin{center}
\vskip-2cm
{\large\bf Pierluigi Colli$^{(1)}$}\\
{\normalsize e-mail: {\tt pierluigi.colli@unipv.it}}\\[.25cm]
{\large\bf Gianni Gilardi$^{(1)}$}\\
{\normalsize e-mail: {\tt gianni.gilardi@unipv.it}}\\[.25cm]
{\large\bf Paolo Podio-Guidugli$^{(2)}$}\\
{\normalsize e-mail: {\tt ppg@uniroma2.it}}\\[.25cm]
{\large\bf J\"urgen Sprekels$^{(3)}$}\\
{\normalsize e-mail: {\tt sprekels@wias-berlin.de}}\\[.45cm]
$^{(1)}$
{\small Dipartimento di Matematica ``F. Casorati'', Universit\`a di Pavia}\\
{\small via Ferrata 1, 27100 Pavia, Italy}\\[.2cm]
$^{(2)}$
{\small Dipartimento di Ingegneria Civile, Universit\`a di Roma ``Tor Vergata''}\\
{\small via del Politecnico 1, 00133 Roma, Italy}\\[.2cm]
$^{(3)}$
{\small Weierstra\ss-Institut f\"ur Angewandte Analysis und Stochastik}\\
{\small Mohrenstra\ss e\ 39, 10117 Berlin, Germany}\\[.8cm]
{\Large{\sl In memory of Enrico Magenes}}\\[.8cm]
\end{center}


\Begin{abstract}
An existence result is proved for a nonlinear diffusion problem of phase-field type, 
consisting of a parabolic system of two partial differential equations, 
complemented by  Neumann homogeneous boundary conditions and initial conditions.  
{\ppg This system is meant to model} 
two-species phase segregation on an atomic lattice 
under the presence of diffusion{\ppg. A similar system} has been recently introduced and analyzed in \cite{CGPS3}{\ppg.  
Both systems conform to the general} theory developed in \cite{Podio}{\ppg: two} parabolic PDEs, interpreted as balances of microforces and microenergy, 
{\ppg are} to be solved {\ppg for  
the} order parameter $\rho$ and the chemical potential~$\mu$.
{\ppg In the system studied in this note, a phase-field equation
in $\rho$   fairly more general than in \cite{CGPS3} is coupled} with a highly nonlinear diffusion equation for $\mu${\ppg,} in which the 
conductivity coefficient {\ppg is allowed to} depend {\ppg nonlinearly} on both variables.
\End{abstract}


\salta

\pagestyle{myheadings}
\newcommand\testopari{\sc Colli \ --- \ Gilardi \ --- \ Podio-Guidugli \ --- \ Sprekels}
\newcommand\testodispari{\sc Existence for a strongly coupled 
Cahn-Hilliard system}
\markboth{\testodispari}{\testopari}

\finqui


\section{Introduction}
\label{Intro}
\setcounter{equation}{0}
In this paper{\ppg,} we prove an existence result 
for the following system in the unknown {\ppg fields} $\mu$ and~$\rho${\ppg:}
\Bsist
  && \bigl( 1 + 2g(\rho) \bigr) \, \dt\mu
  + \mu \, g'(\rho) \, \dt\rho
  - \div \bigl( \kappa(\mu,\rho)\nabla\mu \bigr) = 0{\ppg,}
  \label{Iprima}
  \\
  && \dt\rho - \Delta\rho + f'(\rho) = \mu \, g'(\rho){\ppg,}
  \label{Iseconda}
  \\
  && {\giannipier \big(\kappa(\mu,\rho)\nabla\mu\big)\cdot\nu}|_\Gamma = 0
  \aand
  \dn\rho|_\Gamma = 0{\ppg,}
  \label{Ibc}
  \\
  && \mu(0) = \muz
  \aand
  \rho(0) = \rhoz .
  \label{Icauchy}
\Esist
\Accorpa\Ipbl Iprima Icauchy
{\ppg Each of the} partial differential equations  \accorpa{Iprima}{Iseconda} 
{\giannipier is meant {\piecol to hold} in a three-dimensional bounded domain $\Omega$, endowed with a smooth boundary~$\Gamma$, and in some time interval~$[0,T]$}.
{\giannipier Such a system generalizes the phase-field model of Cahn-Hilliard type} {\ppg studied} 
recently in \cite{CGPS3}{\ppg. Both models are of the type proposed in \cite{Podio}, and aim to describe}
phase segregation of two species (atoms and vacancies, say) on a lattice
in presence of diffusion. 
The state variables are the  {\sl order parameter\/} $\rho$,
interpreted as {\ppg the volume density of one of the two species},  and the {\sl chemical potential\/}~$\mu$. 
For physical reasons, $\mu$~is required to be nonnegative, 
while {\ppg the phase parameter} $\rho$ {\ppg must, as such, obey} 
$0\le \rho\le 1$. {\ppg Here are the features of \cite{CGPS3} that have been generalized. 

Firstly, 
the nonlinearity $f$ considered in \cite{CGPS3} is} a double-well potential defined in~$(0,1)$, 
whose derivative $f'$ {\ppg diverges} at the endpoints $\rho=0$ and $\rho=1${\ppg: e.g., for} 
$f=f_1+f_2$ with $f_2$ smooth, 
one can take  $f_1(\rho)=c\,(\rho\,\log(\rho)+(1-\rho)\,\log(1-\rho))$,  
with $c$ a positive constant. In this paper, we let $f_1$ be a maximal monotone graph from $\erre$ to~$\erre$. {\ppg Consequently, equation~\eqref{Iseconda} }
has to be read as a differential inclusion, 
in which the derivative of the convex part $f_1$ of $f$ is replaced by 
the subdifferential $\beta:=\partial f_1$,~i.e.,
\Beq
  \dt\rho - \Delta\rho + \xi + f_2'(\rho) = \mu g'(\rho)
  \quad \hbox{with} \quad
  \xi \in \beta(\rho);
  \label{Isecondabis}
\Eeq 
{\ppg moreover, since $f_1$ is not required to be smooth, its subdifferential may be multivalued; the selection of $\xi$ in $\beta(\rho)$ is a further difficulty} we face.

{\ppg Secondly, while in \cite{CGPS3} $g(\rho) = \rho$,  here $g$ is any nonnegative-valued smooth function, defined (at least) 
in the domain where $f_1$ and its derivative 
(or rather, its} subdifferential) live. 

{\ppg Thirdly, and this is the most} important novelty{\ppg, conductivity $\kappa$ is not anymore a constant, but rather a positive-valued, continuous, bounded, and possibly nonlinear, function of $\mu$ and $\rho$.  For simplicity, 
we confine ourselves to study the existence of a solution 
under an assumption that guarantees uniform parabolicity, i.e., $\kappa\geq\kmin>0$. We point out that in a recent study \cite{CGPSgen} we let $\kappa$ depend only on $\mu$ and possibly degenerate somewhere.}

{\ppg Finally, relations \eqref{Icauchy} specify the initial conditions for $\mu$ and~$\rho$,
{\giannipier while \eqref{Ibc}} are nothing but homogeneous boundary conditions of Neumann type, involving precisely those boundary operators that match the elliptic differential operators
in \accorpa{Iprima}{Iseconda}. 
}



{\ppg Our} paper is \organiz ed as follows.
In the next section, we state our assumptions and {\ppg our} results.
The existence of a solution is proved in Section~\futuro\ref{Existence}{\ppg ,}
making use of a time-delay approximation and {\ppg of} a number of {\sl a priori estimates\/}, {\ppg that allow}
 us to pass to the limit by compactness and monotonicity techniques. 


\section{Results}
\label{MainResults}
\setcounter{equation}{0}

In this section, we 
describe the mathematical problem under investigation, 
make our assumptions precise, and state our results.
First of all,
we assume $\Omega$ to be a bounded connected open set in $\erre^3$
with smooth boundary~$\Gamma$
({\ppg treating} lower-dimensional cases would require {\ppg only }minor changes).
Next, we fix a final time $T\in(0,+\infty)$
and set:
\Bsist
  && Q := \Omega\times(0,T) , \quad
  \Sigma := \Gamma\times(0,T){\giannipier ,}
  \label{defQS}
  \\
  && V := \Huno,
  \quad H := \Ldue , \quad
  W := \graffe{v\in\Hdue:\ \dn v = 0 \ \hbox{on $\Gamma$}}.
  \label{defspazi}
\Esist
We endow the spaces \eqref{defspazi} with their standard norms,
for which we use a self-explanato\-ry notation 
like $\normaV\cpto${\ppg ; for powers of these spaces, norms are denoted by the same symbols.
We remark that the embeddings $W\subset V\subset H$ are compact,
because $\Omega$ is bounded and smooth. Moreover,} for $p\in[1,+\infty]$, we write $\norma\cpto_p$ for the usual norm
in~$L^p(\Omega)${\ppg; as} no confusion can arise, the symbol $\norma\cpto_p$
is used for the norm in $L^p(Q)$ as well.

\bigskip

First of all, we present the structural assumptions we make.
We require~that{\ppg :}
\Bsist
  \hskip-1.5cm&& \hbox{$\kappa:(m,r)\mapsto\kappa(m,r)$ is continuous from $[0,+\infty)\times\erre$ to $\erre$}{\ppg ,}
  \label{hpk}
  \\
  \hskip-1.5cm&& \hbox{the partial derivatives $\partial_r\kappa$ and $\partial_r^2\kappa$
  exist and are continuous}{\ppg ,}
  \label{hpdk}
  \\
  \hskip-1.5cm&& \kmin,\kmax \in (0,+\infty){\giannipier ,}
  \label{hpcost}
  \\
  \hskip-1.5cm&& {\giannipier \kmin \leq \kappa(m,r) \leq \kmax, \enskip
  |\partial_r\kappa(m,r)| \leq \kmax,
  \enskip 
  |\partial_r^2\kappa(m,r)| \leq \kmax \quad 
  \hbox{for $m\geq0$ and $r\in\erre$}}
  \label{hpkbis}
  \\
  \hskip-1.5cm&& 
  {\pier f = f_1 + f_2 \,, \quad f_1:\erre \to [0,+\infty], \quad f_2 :\erre \to \erre,} \quad 
  g:\erre \to [0,+\infty){\giannipier ,}
  \label{hpfg}
  \\
  \hskip-1.5cm&& \hbox{$f_1$ is convex, proper, l.s.c. and $f_2$ and $g$ are $C^2$ functions}{\giannipier ,}
  \label{hpfuno}
  \\
  \hskip-1.5cm&& \hbox{$f_2'$, $g$, and $g'$ are Lipschitz continuous}.
  \label{hpfdueg}
\Esist
{\ppg For convenience, we set:}
\Bsist
  \hskip-1.5cm&& \kappa' := \partial_r \kappa , \quad
  \kappa'' := \partial_r^2 \kappa , \quad
  \beta := \partial f_1 \,,
  \aand
  \pi := f_2' {\ppg ;}
  \label{defbp}
  \\
  \hskip-1.5cm&& K(m,r) := \int_0^m \!\!\! \kappa(s,r) \ds, \
  \Ku(m,r) := \int_0^m \!\!\! \kappa'(s,r) \ds, \
  \Kd(m,r) := \int_0^m \!\!\! \kappa''(s,r) \ds 
  \non
  \\
  \hskip-1.5cm&& \quad \hbox{for $m\geq0$ and $r\in\erre$}{\ppg ;}
  \label{defK}
\Esist
\Accorpa\Hpstruttura hpk defK
and {\ppg we} write $D(f_1)$ and $D(\beta)$ for the effective domains
of~$f_1$ and~$\beta$, respectively.
Clearly{\ppg , thanks to \eqref{hpkbis},}
\Beq
  \max \{ |K(m,r)| , |\Ku(m,r)| , |\Kd(m,r)| \}
  \leq \kmax m
  \quad \hbox{for every $m\geq0$ and $r\in\erre $}.
  \label{lingr}
\Eeq
We also note that the structural assumption{\ppg s} of~\cite{CGPSgen}
are fulfilled if $\kappa$ only depends on~$m${\ppg ,} 
and that a strong singularity in equations \eqref{Iseconda} for~$\rho$
is allowed. {\ppg At variance with ~\cite{CGPSgen},} equation \eqref{Iprima} for $\mu$ is {\ppg here} uniformly parabolic,
since $g$ is nonnegative and $\kappa$ is bounded away from zero.

{\pier
\Brem
\label{justif}
Note that any convex, proper, l.s.c.\ {\giannipier function} is bounded from below by an affine function 
(cf., e.g., \cite[Prop.~2.1, p.~51]{Barbu}){\ppg ,}
so that the assumption $f_1 \geq 0 $ looks reasonable{\ppg , because} one can suitably modify the smooth perturbation~$f_2$. 
Moreover, let us point out that the other positivity condition{\ppg ,} $g\geq 0${\ppg ,} 
is just needed on the set~$D(\beta)$, 
{\ppg while $g$ can be extended }outside of $D(\beta)$ accordingly. 
\Erem
}

As {\ppg to} initial data,
we~require~that:
\Bsist
  && \muz \in V \cap \Linfty , \quad
  \rhoz \in W , \quad
  \muz \geq 0
  \aand \rhoz \in D(\beta) \quad \aeO{\ppg ;}
  \label{hpzero}
  \\
  && \hbox{there exists $\xi_0\in H$ such that $\xi_0\in\beta(\rhoz)$ \aeO}.
  \label{hprhozbis}
\Esist
\Accorpa\Hpdati hpzero hprhozbis
{\ppg Since} $f_1$ is convex and $f_2$ smooth,
the above assumptions imply {\ppg that} $f(\rhoz)\in\Luno$.

{\ppg As to the a priori regularity we require 
for any solution $(\mu,\rho,\xi)$, we begin to observe} that, for any given~$\mu$, equation~\eqref{Isecondabis}
has the form of a standard phase{\ppg -}field equation.
Therefore, it is natural to look for pairs $(\rho,\xi)$ that~satisfy
\Bsist
  && \rho \in \W{1,\infty}H \cap \H1V \cap \L\infty W
  \label{regrho}{\ppg ,}
  \\
  && 
  {\giannipier \xi \in \L\infty H  .}
  \label{regxi}
\Esist
{\giannipier Note} {\ppg that the Neumann boundary condition 
{\giannipier for $\rho$} 
has been incorporated into \eqref{regrho}} 
{\giannipier (cf.~$\eqref{defspazi}_3$).} 
{\ppg Next, as to} $\mu$,
we require~that
\Bsist
  \hskip-1.5cm&& \mu \in \H1H \cap \L\infty V \cap \LQ\infty , \quad
  \mu \geq 0 \quad \aeQ{\ppg ,}
  \label{regmu}
  \\
  \hskip-1.5cm&& \div \bigl( \kamurho\nabla\mu \bigr) \in \L2H 
  \aand
  {\giannipier \bigl( \kamurho\nabla\mu \bigr) \cdot \nu = 0
  \quad \aeS,}
  \label{regdiv}
\Esist
\Accorpa\Regsoluz regrho regdiv
and note that we can expect that $\mu\in\L2W$
(from the regularity theory of elliptic equations)
only if the function $\kappa$ is smooth with respect to both variables.
Nevertheless, \eqref{regmu} and the regularity of the divergence
are sufficient to write the Neumann boundary condition as done in~\eqref{regdiv}.
We also observe that 
\Beq
  \rho \in \C0{\Cx0} = \CQ{\ppg ,}
  \label{rhocont}
\Eeq
{\ppg as a direct consequence of} \eqref{regrho} and the compact embedding $W\subset\Cx0$
(see, e.g., \cite[Sect.~8, Cor.~4]{Simon}),
whence $g'(\rho)\in \CQ$.
Thus, under all {\ppg of} the above requirements,
we can write the system of equations and the initial condition in the following strong form
\Bsist
  & \bigl( \coeff \bigr) \dt\mu + \mu \, g'(\rho) \, \dt\rho
  - \div \bigl( \kamurho \nabla\mu \bigr) = 0
  & \quad \aeQ{\ppg ,}
  \label{prima}
  \\
  & \dt\rho - \Delta\rho + \xi + \pi(\rho)
  = \mu \, g'(\rho)
  \aand \xi \in \beta(\rho)
  & \quad \aeQ{\ppg ,}
  \label{seconda}
  \\
  & \mu(0) = \muz
  \aand
  \rho(0) = \rhoz
  & \quad \aeO{\ppg .}
  \label{cauchy}
\Esist
\Accorpa\Pbl prima cauchy

Here is our existence result.

\Bthm
\label{Esistenza}
Assume {\ppg that both} \Hpstruttura\ and \Hpdati\ {\ppg hold}.
Then, there exists at least a triplet $(\mu,\rho,\xi)$ satisfying \Regsoluz\ and solving
problem~\Pbl.
\Ethm

{\ppg This is the only result we prove} in the present paper. {\ppg We note, however, that}
 the uniqueness {\ppg result} obtained in~\cite{CGPSgen}
still holds here{\ppg ,} provided that $\kappa$ is {\ppg taken} constant and $\muz$  smoother.
For the reader's convenience,
we \summariz e the results of~\cite{CGPSgen}.

\Bthm
\label{Precedenti}
Assume {\ppg that both} \Hpstruttura\ and \Hpdati\ {\ppg hold}, {\ppg and, moreover, that} $\muz\in W$ and $\kappa=1$.
Then, {\ppg there is a unique} triplet $(\mu,\rho,\xi)$ 
satisfying \Regsoluz\ and solving
problem~\Pbl\, and its component $\mu$ 
enjoys the {\ppg following} regularity property{\ppg :}
\Beq
  \mu \in \W{1,p}H \cap \L pW
  \quad \hbox{for every $p\in[1,+\infty)$}.
  \label{piuregmu}
\Eeq
\Ethm

Throughout {\ppg the} paper,
we {\ppg make use of} some \wk\ embeddings
of Sobolev type,
namely{\ppg ,}
$V\subset\Lx p$ for $p\in[1,6]$,
{\giannipier together} with the related Sobolev inequality
\Beq
  \norma v_p \leq C \normaV v
  \quad \hbox{for every $v\in V$ and $1\leq p \leq 6$}{\ppg ,}
  \label{sobolev}
\Eeq
and $\Wx{1,p}\subset\Cx0$ for $p>3${\ppg , together} with 
\Beq
  \norma v_\infty \leq C_p \norma v_{\Wx{1,p}}
  \quad \hbox{for every $v\in\Wx{1,p}$ and $p>3$}.
  \label{sobolevbis}
\Eeq
In \eqref{sobolev}, $C$ depends {\ppg only} on~$\Omega$, 
while $C_p$ in \eqref{sobolevbis} depends {\ppg also} on~$p$.
In particular, the continuous embedding 
$W\subset\Wx{1,6}\subset\Cx0$ {\pier holds}.
Some of the previous embeddings are in fact compact.
This is the case for $V\subset\Lq$ and $W\subset\Cx0$.
We also account for the corresponding inequality
\Beq
  \norma v_4 \leq \eps \normaH{\nabla v} + C_\eps \normaH v
  \quad \hbox{for every $v\in V$ and $\eps>0$}
  \label{compact}
\Eeq
where $C_\eps$ depends on $\Omega$ and~$\eps$, only.
Furthermore, we repeatedly make use of the notation
\Beq
  Q_t := \Omega \times (0,t)
  \quad \hbox{for $t\in[0,T]$}{\ppg ,}
  \label{defQt}
\Eeq
of the \wk\ \holder\ inequality{\ppg , and of} the elementary Young inequality
\Beq
  ab \leq \eps a^2 + \frac 1{4\eps} \, b^2
  \quad \hbox{for every $a,b\geq 0$ and $\eps>0$}.
  \label{young}
\Eeq
Finally, {\ppg again} throughout the paper,
we use a small-case italic $c$ for different constants, that
may only depend 
on~$\Omega$, the final time~$T$, the shape of the nonlinearities $f$ and~$g$, 
and the properties of the data involved in the statements at hand; 
a~notation like~$c_\eps$ signals a constant that depends also on the parameter~$\eps$. 
The reader should keep in mind that the meaning of $c$ and $c_\eps$ might
change from line to line and even in the same chain of inequalities, 
whereas those constants we need to refer to are always denoted by 
capital letters, just like $C$ in~\eqref{sobolev}.


\section{Existence}
\label{Existence}
\setcounter{equation}{0}

In this section, we prove Theorem~\ref{Esistenza}{\ppg ,} which ensures the existence of a solution.
{\pier Although} our proof follows the argument {\ppg in}~\cite{CGPS3} and~\cite{CGPSgen} closely,
we present the whole argument{\ppg ,} and sometimes give some details,
since the changes with respect to the quoted papers are spread {\ppg over the whole exposition}.
{\ppg Our} starting point is an approximating problem{\ppg ,}
which is still based on {\ppg introducing} a time delay in the \rhs\
of~\eqref{seconda}.
{\ppg Precisely}, we define the translation operator $\T:\L1H\to\L1H$
depending
on a time step $\tau>0$ by setting, for $v\in\L1H$ and \aat,
\Beq
  (\T v)(t) := v(t-\tau)
  \quad \hbox{if $t>\tau$}
  \aand
  (\T v)(t) := \muz
  \quad \hbox{if $t<\tau$}
  \label{defT}
\Eeq 
(the same notation $\T v$ will be used {\ppg also} for a function $v$
that is defined in some subinterval $[0,T']$ of~$[0,T]$){\ppg . At bottom, what we do is to} replace $\mu$ by $\T\mu$ in~\eqref{seconda}.
However, since it is not obvious that we can keep {\ppg $\mu$ positive},
we extend $\kappa$ to a function $\kappat:\erre\to\erre$
satisfying similar properties.
Moreover, we assume that the analogue of \eqref{hpkbis} holds 
for $\kappat$ and its derivatives{\ppg ,} with the same constants $\kmin$ and~$\kmax$
(we~replace $\kmin$ and $\kmax$ by $2\kmin$ and $\kmax/2$ in the original~\eqref{hpkbis}
if necessary).
So, the approximating problem consists of the equations
\Bsist
  & \bigl( \coefftau \bigr) \, \dt\mut
  + \mut \, g'(\rhot) \, \dt\rhot
  - \div \bigl( \kappat(\mut,\rhot)\nabla\mut \bigr) = 0
  & \quad \aeQ{\ppg ,}
  \label{primatau}
  \\
  & \dt\rhot - \Delta\rhot + \xit + \pi(\rhot)
  = (\T\mut) \, g'(\rhot)
  \aand \xit \in \beta(\rhot)
  & \quad \aeQ{\ppg ,}
  \label{secondatau}
\Esist
complemented {\ppg by} the initial and boundary conditions
\Beq
  \mut(0) = \muz \,,\quad
  \rhot(0) = \rhoz \,,\quad
  \dn\rhot|_\Sigma = 0 ,\quad
  (\kappat(\mut,\rhot)\nabla\mut\cdot\nu|_\Sigma = 0.
  \label{caubctau}
\Eeq
For convenience, we allow $\tau$ to take just discrete values,
namely, $\tau=T/N$, where $N$ is any positive integer.

\Blem
\label{EsistPbltau}
The approximating problem has a solution $(\mut,\rhot,\xit)$ satisfying
the analogue of \Regsoluz.
\Elem

\Bdim
We confine ourselves to give a sketch.
As in~\cite{CGPS3}, we inductively solve $N$ problems
on the time intervals $I_n=[0,t_n]:=[0,n\tau]$, $n=1,\dots,N$,
by constructing the solution directly on the whole of~$I_n$ at each step.
Namely, given $\munmu$, which is defined in $\Omega\times I_{n-1}$,
we note that $\T\munmu$ is well defined and known 
in $\Omega\times I_n$
and solve the boundary value problem for $\rhon$
given by the phase field equations
\Beq
  \dt\rhon - \Delta\rhon + \xin + \pi(\rhon)
  = (\T\mun) \, g'(\rhon)
  \aand \xin \in \beta(\rhon)
  \quad \hbox{in $\Omega\times I_n$}
  \label{secondan}
\Eeq
complemented {\ppg by} the boundary and initial conditions just mentioned for~$\rhot$.
Such a problem is quite standard and has a unique solution $\rhon$
in a proper functional framework.
Now, we observe that the function
\Beq
  \hat\kappa: (x,t,m) \mapsto \kappat(m,\rhon(x,t)) , \quad
  (x,t) \in \Omega \times I_n, \quad m \in \erre
  \non
\Eeq
is a Carath\'eodory function satisfying 
$\kmin\leq\hat\kappa\leq\kmax$ on its domain,
so that the equation
\Beq
  \bigl( \coeffn \bigr) \, \dt\mun
  + \mun \, g'(\rhon) \, \dt\rhon
  - \div \bigl( \kappat(\mun,\rhon)\nabla\mun \bigr) = 0
  \quad \hbox{in $\Omega\times I_n$}
  \label{priman}
\Eeq
in the unknown function $\mun$ is uniformly parabolic 
({\pier let also recall} that $g$ is nonnegative).
Thus, we can solve the problem obtained by complementing \eqref{priman}
with the boundary and initial conditions prescribed for~$\mut$.
Therefore, the problem to be solved has a unique solution in a proper space{\ppg ,}
provided that the coefficient $g'(\rhon)\dt\rhon$ is not too irregular.
So, we should prove that, step by step, we get the right regularity
for $\rhon$ and~$\mun$.
This could be done by induction, as in~\cite{CGPS3}, with some modifications
due to our more general framework.
We omit this detail and just observe that the needed a~priori estimates 
are close (and even simpler{\ppg ,} since $\tau$ is fixed here)
to~the ones we perform later on in order to let $\tau$ go to zero.
The final point is $\mun\geq0$.
We give the proof in detail.
We multiply equation~\eqref{priman} by $-\mun^-:=-(-\mun)^+$,
the negative part of~$\mun$,
and integrate over~$Q_t$ with any $t\in I_n$.
We observe that
\Beq
  \bigl[ \bigl( \coeffnt \bigr) \, \dt\mun + \mun \, g(\rhon) \, \dt\rhon \bigr] \, (-\mun^-)
  = \frac 12 \, \dt \bigl( (\coeffn) \, |\mun^-|^2 \bigr).
  \non
\Eeq
Hence, by using $\muz\geq0$ 
and owing to the boundary condition, we have
\Bsist
  \frac 12 \iO (\coeffnt) \, |\mun^-(t)|^2
  + \intQt \kappat(\mun,\rhon) |\nabla\mun^-|^2 
  = 0 
  \quad \hbox{for every $t\in I_n$}.
  \non
\Esist
As {\ppg both} $g$ and $\kappat$ are nonnegative,
this implies $\mun^-=0$, that is, $\mun\geq0$ a.e.\ in~$\Omega\times I_n$.
Once all this is checked, the finite sequence
$(\mun,\rhon,\xin)$, $n=1,\dots,N$, is constructed 
and it is clear that 
a solution to the approximating problem we are looking for
is simply obained by taking $n=N$.
\Edim

{\pier Thus, we fix} a solution $(\mut,\rhot,\xit)$ for each~$\tau$.
We note that{\ppg , a~posteriori,} we can replace $\kappat$ by $\kappa$ in \eqref{primatau}, 
since the component $\mut$ of our solution is nonnegative.
Our aim is to let $\tau$ go to zero{\ppg , so as} to obtain
a solution as stated in Theorem~\ref{Esistenza}.
Our proof uses compactness arguments and thus relies on a number of 
uniform (with respect to~$\tau$) a~priori estimates.
In order to make the formulas to come more readable, 
we write $\mu$ and $\rho$ rather than $\mut$ and $\rhot$ in the calculations.

\Brem
\label{Noformal}
{\ppg Sometimes, when} deriving our a priori estimates, we proceed formally.
However, our procedure{\ppg s} can be made rigorous.
For instance, one can get more regularity for the approximating problem 
by \regulariz ing~$\kappat$ and the initial data, if necessary.
Moreover, {\ppg local regularization} is often sufficient.
{\ppg Consider, e.g.,} equation~\eqref{secondatau} and {\pier rewrite} it in the form
\Beq
  -\Delta\rho + \rho + \beta(\rho)
  \ni h := \rho - \dt\rho - \pi(\rho) + (\T\mu)g'(\rho) .
  \label{noformal}
\Eeq
Now,  (here $t$ is just a parameter) the elliptic equation:
\Beq
  -\Delta u + u + w = h
  \aand
  w \in \beta(u){\ppg ,}
  \non
\Eeq
complemented {\ppg by} homogeneous Neumann boundary condition{\ppg ,}
yields a well-posed problem{\ppg ;} as is well known,
its solution $(u,w)$ is the limit of the much smoother pair $(u_\eps,\beta_\eps(u_\eps))$,
where $\beta_\eps$ is a \regulariz ation of $\beta$ of Yosida type
(see, e.g., \cite[p.~28]{Brezis}; 
see also the proof of Lemma~3.1 of~\cite{CGPS3} for a further \regulariz ation)
and~$u_\eps$ is the solution of the analogous boundary value problem~for
\Beq
  -\Delta u_\eps + u_\eps + \beta_\eps(u_\eps)
  = h .
  \non
\Eeq
On the other hand, we have $(u,w)=(\rho,\xi)$ by~\eqref{noformal}.
Therefore, it is essentially correct to {\pier regard} $\beta$ 
as {\ppg if} it were a smooth function in the original equation \eqref{secondatau}{\ppg ,}
and treat such equation like a more regular one
(e.g., by differentiating it or taking irregular functions as test functions).
\Erem

\step 
First a priori estimate

We test \eqref{primatau} by $\mut$ and observe that
\Beq
  \bigl[ \bigl( \coeff \bigr) \, \dt\mu
  + \mu \, g'(\rho) \, \dt\rho \bigr] \mu
  = \frac 12 \, \dt \bigl[ (\coeff) \mu^2 \bigr].
  \non
\Eeq
Thus, by integrating over~$(0,t)$, where $t\in[0,T]$ is arbitrary,
we obtain
\Beq
  \iO \bigl( \coefft \bigr) |\mu(t)|^2 + \intQt \kamurho |\nabla\mu|^2
  = \iO (1+2g(\rhoz)) \muz^2 \,.
  \non
\Eeq
Hence, we recall that $g\geq0$ and $\kappat\geq\kmin>0${\ppg ,}
and conclude that
\Beq
  \norma\mut_{\L\infty H} + \norma{\mut}_{\L2V} \leq c .
  \label{primastima}
\Eeq

\step Consequence

The Sobolev inequality~\eqref{sobolev}, estimate \eqref{primastima},
and \eqref{lingr}{\ppg ,} imply that
\Beq
  \norma\mut_{\L2{\Lx6}} 
  + \norma{\psi(\mut,\rhot)}_{\L\infty H\cap\L2{\Lx6}} \leq c
  \quad \hbox{with $\psi=K,\Ku,\Kd$}.
  \label{daprima}
\Eeq
{\piecol Another implication of \eqref{primastima}, along with \eqref{defT} and \eqref{hpzero}, is}
\Beq
 {\piecol  \norma{\T \mut}_{\L\infty H} + \norma{\T \mut}_{\L2V} \leq c .}
  \label{pristitau} 
\Eeq

\step 
Second a priori estimate

We {\pier add $\rhot$ to both sides of \eqref{secondatau} and test by} $\dt\rhot$.
We obtain{\ppg :}
\Bsist
  && \intQt |\dt\rho|^2
  + {\pier \frac 12 \, \norma{\rho(t)}_V^2 }
  + \iO f_1(\rho(t))
  \non
  \\
  && = \frac 12 \iO |\nabla\rhoz|^2
  + \iO f(\rhoz)
  + {\pier \frac 12 \iO \tonde{\rho^2 (t) + 2 f_2(\rho(t))} }
  + \intQt g'(\rho) (\T\mu) \dt\rho 
  \non
  \\
  && \leq c + {\pier c \iO |\rho (t)|^2} + \frac 14 \intQt |\dt\rho|^2
  + c \norma{\T\mu}_{\L\infty H}^2 ,
  \non
\Esist
{\ppg  for every $t\in[0,T]$.}
{\pier Thanks to the chain rule and {\gianni the Young inequality} \eqref{young}, we have{\ppg :}
\Beq
  c \iO |\rho(t)|^2 
  \leq c \iO |\rhoz|^2 + \frac 14  \intQt |\dt\rho|^2  
  + c\int_0^t \norma{\rho(s)}_H^2\, ds.
  \non
\Eeq
Then, as $f_1$ is nonnegative, by accounting for \eqref{primastima}, with the help of the Gronwall lemma we infer that
\Beq
  \intQt |\dt\rho|^2
  +  \norma{\rho(t)}_V^2
  + \iO  f_1(\rho(t)) \leq c.
\non
\Eeq
Thus, we conclude~that}
\Beq
  \norma\rhot_{\H1H\cap\L\infty V} \leq c
  \aand
  \norma{f(\rhot)}_{\L\infty\Luno} \leq c.
  \label{secondastima}
\Eeq

\step 
Third a priori estimate

We proceed formally (see Remark~\ref{Noformal}).
We {\pier rewrite} \eqref{secondatau} as
\Beq
  - \Delta\rho + \beta(\rho) 
  = h := - \dt\rho - \pi(\rho) + (\T\mu)g'(\rho){\ppg ,}
  \label{perterza}
\Eeq
and multiply {\ppg this relation} by either $-\Delta\rho$ or $\beta(\rho)$.
By doing that, we derive an estimate for both $\Delta\rho$ and $\beta(\rho)$
and {\ppg we} can use the regularity theory for elliptic equations.
We conclude~that
\Beq
  \norma\rhot_{\L2W} \leq c
  \aand
  \norma\xit_{\L2H} \leq c.
  \label{terzastima}
\Eeq

\step 
Fourth a priori estimate

Our aim is to improve estimates~\eqref{secondastima} and~\eqref{terzastima}.
To do that, we proceed formally, at least at the beginning {\ppg (our  procedure could be made completely rigorous,}
as sketched in Remark~\ref{Noformal}).
We start from {\ppg an} estimate {\ppg coming from} the theory of maximal monotone operators {\ppg \cite{Brezis}, namely,}
\Beq
  \normaH{\dt u(0)}
  \leq \normaH{\psi(0)+\Delta\rhoz}
  + \min_{\eta\in\beta(\rhoz)} \normaH\eta{\ppg ,}
  \label{perquartaZ}
\Eeq
for the unique solution $(u,\omega)$ to the equations {\pier (cf.~\eqref{secondatau})}
\Beq
  \dt u - \Delta u + \omega = \psi := g'(\rho) \T\mu - \pi(\rho)
  \aand
  \omega \in \beta(u){\ppg ,}
  \non
\Eeq
complemented {\ppg by} the same initial and boundary conditions prescribed for~$\rho$.
{\pier Note that in \eqref{perquartaZ} $\beta$ is understood as the induced maximal 
monotone operator from $H$ to $H$. By observing} 
that $(u,\omega)=(\rho,\xi)$, applying \eqref{perquartaZ},
and combining with our assumptions on~$\rhoz$
(see \eqref{hprhozbis}, in particular),
we obtain{\ppg :}
\Beq
  \normaH{\dt\rhot(0)}
  \leq c \bigl( \normaH\muz + \normaW\rhoz + 1 + \normaH{\xi_0} \bigl)
  = c .
  \label{finequartaZ}
\Eeq
We use \eqref{finequartaZ} in the calculation we are {\ppg about to start: once again, we proceed}
 formally, 
and {\ppg write} $\xi=\beta(\rho)$ 
as {\ppg if} $\beta$ were a smooth function.
We differentiate \eqref{secondatau} with respect to time
and test the equation {\pier obtained for} $\dt\rho$.
We~{\ppg find:}
\Bsist
  && \frac 12 \iO |\dt\rho(t)|^2
  + \intQt |\nabla\dt\rho|^2
  + \intQt \beta'(\rho) |\dt\rho|^2
  \non
  \\
  && = \frac 12 \iO |(\dt\rho)(0)|^2 
  - \intQt \pi'(\rho) |\dt\rho|^2
  + \intQt g''(\rho) {\piecol(\T\mu)} |\dt\rho|^2
  \qquad
  \non
  \\
  &&\quad {} 
  + \intQt g'(\rho) \dt(\T\mu) \, \dt\rho
  \non
  \\
  && \leq \frac 12 \iO |(\dt\rho)(0)|^2 
  + c \intQt (1+{\piecol \T \mu}) |\dt\rho|^2
  + \intQt g'(\rho) \dt(\T\mu) \, \dt\rho .
  \label{perquarta}
\Esist
We treat each term on the \rhs, separately.
The first one is estimated by~\eqref{finequartaZ}.
In order to deal with the second one,
we account for the \holder\ inequality,
{\piecol \eqref{pristitau}}, 
the compact embedding $V\subset\Lq$ (see \eqref{compact}),
and \eqref{secondastima}. 
We~obtain{\ppg :}
\Bsist
  && \intQt (1+{\piecol \T \mu}) |\dt\rho|^2
  \leq \iot \normaH{1+{\piecol (\T \mu )}(s)} \norma{\dt\rho(s)}_4^2 \ds
  \non
  \\
  && \leq c \iot \norma{\dt\rho(s)}_4^2 \ds
  \leq \eps \intQt |\nabla\dt\rho|^2
  + c_\eps \intQt |\dt\rho|^2
  \qquad
  \non
  \\
  && \leq \eps \intQt |\nabla\dt\rho|^2
  + c_\eps \,,
  \label{finequartaA}
\Esist
{\ppg for every~$\eps>0$}. Let us come to the last term of~\eqref{perquarta}.
{\ppg Firstly, on recalling} that $\dt\T\mu=0$ in $(0,\tau)$ by the definition of~$\T${\ppg , we}
compute $\dt\mu$ from~\eqref{primatau}.
Then{\ppg ,} we integrate by parts and {\ppg have repeated recourse} to \holder, Sobolev, and Young 
inequalities. We {\pier deduce that}
\Bsist
  \hskip-1cm&& \intQt g'(\rho) \dt(\T\mu) \, \dt\rho
  = \intQtmt \dt\mu(s) \, \dt g(\rho(s+\tau)) \ds
  \non
  \\
  \hskip-1cm&& = \intQtmt \frac 1{\coeffs} \,
  \bigl[ \div \bigl(\kamurhos\nabla\mu(s)\bigr) - \mu(s) g'(\rho(s)) \dt\rho(s) \bigr]
  \dt g(\rho(s+\tau)) \ds
  \non
  \\
  \hskip-1cm&& = \intQtmt \kamurhos \nabla\mu(s) \cdot \nabla \,  \frac 
  {\dt g(\rho(s+\tau))} {\coeffs} \ds
  \non
  \\
 \hskip-1cm && \qquad {}
  - \intQtmt \frac {g'(\rho(s)) g'(\rho(s+\tau)) }{\coeffs} \, \mu(s) \dt\rho(s) \dt\rho(s+\tau) \ds {\ppg \,.}
  \label{perquartaB}
\Esist
{\ppg We }treat the last two integrals separately,
by using our structural assumptions.

{\ppg As to the former, we have:} 
\Bsist
  && \intQtmt \kamurhos \nabla\mu(s) \cdot \nabla \, 
  \frac{\dt g(\rho(s+\tau))}{\coeffs} \ds
  \non
  \\
  && = {\piecol \intQtmt \kamurhos \nabla\mu(s) \cdot \nabla \, 
  \frac{g'(\rho(s+\tau)) \dt \rho(s+\tau)}{\coeffs} \ds }
  \non
  \\  
  && \leq c \intQtmt |\nabla\mu(s)| \, |\nabla\dt\rho(s+\tau)| \ds
  \non
  \\
  && \quad {}
  + c \intQtmt |\nabla\mu(s)| \, |\nabla\rho(s )| 
  \, |\dt\rho(s+\tau)| \ds
  \non
  \\
  && \quad 
  {\piecol {} + c \intQtmt |\nabla\mu(s)| \, |\nabla \rho(s+\tau)| 
  \, |\dt\rho(s+\tau)| \ds} .\label{perquartaBA}
\Esist
{\ppg Moreover, thanks to~\eqref{primastima}, we {\piecol infer}:}
\Bsist
  && \intQtmt |\nabla\mu(s)| \, |\nabla\dt\rho(s+\tau)| \ds
  \leq \eps \intQt |\nabla\dt\rho|^2
  + c_\eps \intQt |\nabla\mu|^2
  \non
  \\
  && \leq \eps \intQt |\nabla\dt\rho|^2 + c_\eps\,,
  \label{perquartaBAA}
\Esist
 for every $\eps\in(0,1)$.
On the other hand{\ppg ,} {\piecol we also have:}
\Bsist
  && \intQtmt |\nabla\mu(s)| \, |\nabla\rho(s)| \, |\dt\rho(s+\tau)| \ds
  \non
  \\
  && \leq \iotmt \norma{\nabla\mu(s)}_2 \norma{\nabla\rho(s)}_4 \norma{\dt\rho(s+\tau)}_4 \ds
  \non
  \\
  && \leq \eps \iot \normaV{\dt\rho(s)}^2 \ds
  + c_\eps {\piecol\iotmt} \normaH{\nabla\mu(s)}^2 \normaV{\nabla\rho(s)}^2 \ds
  \non
  \\
  && \leq \eps \intQt |\nabla\dt\rho|^2
  + c \intQt |\dt\rho|^2
  + c_\eps {\piecol\iotmt} \normaH{\nabla\mu(s)}^2 \normaV{\nabla\rho(s)}^2 \ds
  \non
  \\
  && \leq \eps \intQt |\nabla\dt\rho|^2 + c 
  + c_\eps {\piecol\iotmt} \normaH{\nabla\mu(s)}^2 \normaV{\nabla\rho(s)}^2 \ds
\Esist
{\ppg (in the last inequality, \eqref{secondastima} has been used)}.
Now, we improve the estimate just obtained
by owing to the regularity theory for linear elliptic equations{\ppg, as well as}
to estimates \eqref{primastima} and \eqref{secondastima}.
For {\ppg each} fixed~$\piecol{ s \in (0,T)}$, we have
\Bsist
  && \normaV{\nabla\rho(s)}^2
  \leq c \bigl( \normaV{\rho(s)}^2 + \normaH{\Delta\rho(s)}^2 \bigr)
  \non
  \\
  && \leq c + c \normaH{-\dt\rho(s) - \pi(\rho(s)) + g'(\rho(s)) \T\mu(s)}^2
  \leq \normaH{\dt\rho(s)}^2 + c.
  \non
\Esist
Therefore, the above estimate becomes
\Bsist
  && \intQtmt |\nabla\mu(s)| \, |\nabla\rho(s)| \, |\dt\rho(s+\tau)| \ds
  \non
  \\
  && \leq \eps \intQt |\nabla\dt\rho|^2 
  + c_\eps \iot \norma{\nabla\mu(s)}_2^2 \, \normaH{\dt\rho(s)}^2 \ds
  + c_\eps \,.
  \label{perquartaBAB}
\Esist
{\piecol Analogously, one shows that
\Bsist
  && \intQtmt |\nabla\mu(s)| \, |\nabla\rho(s+\tau)| \, |\dt\rho(s+\tau)| \ds
  \non
  \\
  && \leq \eps \intQt |\nabla\dt\rho|^2 
  + c_\eps \iot \norma{\nabla(\T \mu)(s)}_2^2 \, \normaH{\dt\rho(s)}^2 \ds
  + c_\eps \,. 
  \label{perquartainpiu}
\Esist
}
Thus, by collecting \eqref{perquartaBAA} and {\piecol \eqref{perquartaBAB}--\eqref{perquartainpiu}},
we deduce that \eqref{perquartaBA} yields{\ppg :}
\Bsist
  \hskip-1cm&& \intQtmt \kamurhos \nabla\mu(s) \cdot \nabla \, \frac {\dt g(\rho(s+\tau))} {\coeffs} \ds
  \non
  \\
  \hskip-1cm&& \leq \eps \intQt |\nabla\dt\rho|^2 
  + c_\eps \iot {\piecol \Big( \normaH{\nabla\mu(s)}^2 
  + \normaH{\nabla(\T \mu)(s)}^2 \Big)}
  \normaH{\dt\rho(s)}^2 \ds
  + c_\eps \,,
  \label{finequartaBA}
\Esist
{\ppg for every $\eps>0$. 

We now take up the last integral in}~\eqref{perquartaB}.
By using the compactness inequality \eqref{compact} and \eqref{primastima},
we have{\ppg :}
\Bsist
  && - \intQtmt \frac {g'(\rho(s)) g'(\rho(s+\tau))}{\coeffs} \, \mu(s) \dt\rho(s) \dt\rho(s+\tau) \ds 
  \non
  \\
  && \leq c \iotmt \norma{\mu(s)}_4 \norma{\dt\rho(s+\tau)}_4 \norma{\dt\rho(s)}_2 \ds
  \non
  \\
  && \leq \eps \iotmt\normaV{\dt\rho(s+\tau)}^2 \ds 
  + c_\eps \iot \norma{\mu(s)}_4^2 \normaH{\dt\rho(s)}^2 \ds
  \non
  \\
  && \leq  {\pier \eps \intQt |\nabla\dt\rho|^2 } + c
  + c_\eps \iot \norma{\mu(s)}_4^2 \normaH{\dt\rho(s)}^2 \ds .
  \label{finequartaBB}
\Esist
Therefore, due to \eqref{finequartaBA} and \eqref{finequartaBB},
\eqref{perquartaB} becomes{\ppg :}
\Bsist
   &&\intQt g'(\rho) \dt(\T\mu) \, \dt\rho
  \leq {\pier 2 }\eps \intQt |\nabla\dt\rho|^2 
  + c_\eps \iot \normaH{\nabla\mu(s)}^2 \normaH{\dt\rho(s)}^2 \ds
  \non
  \\
  && \hskip4cm {}
  + c_\eps \iot \norma{\mu(s)}_4^2 \normaH{\dt\rho(s)}^2 \ds 
  + c_\eps \,.
  \label{finequartaB}
\Esist
At this point, we combine \eqref{finequartaZ}, \eqref{finequartaA}, {\ppg and} \eqref{finequartaB}
with \eqref{perquarta}{\ppg ,} and {\ppg we} choose $\eps$ small enough.
{\ppg Since} the last integral on the \lhs\ of \eqref{perquarta} is nonnegative{\ppg , because} $f_1$ is convex,
we obtain{\ppg :}
\Bsist
  && \iO |\dt\rho(t)|^2
  + \intQt |\nabla\dt\rho|^2
  \leq c \iot \phi(s) \normaH{\dt\rho(s)}^2 \ds
  + c{\ppg ,}
  \non
  \\
  && \quad \hbox{where} \quad
  \phi(s) := \normaH{\nabla\mu(s)}^2 
{\piecol {}+ \normaH{\nabla(\T \mu)(s)}^2 }
+ \norma{\mu(s)}_4^2 \,.
  \non
\Esist
As $\phi\in L^1(0,T)$ by {\piecol \eqref{primastima}--\eqref{pristitau}},
we can apply the Gronwall lemma and conclude~that
\Beq
  \norma{\dt\rhot}_{\L\infty H\cap\L2V} \leq c.
  \label{quartastima}
\Eeq

\step
Consequence

{\giannipier Note} that
$-\Delta\rhot+\xit=-\dt\rhot+g'(\rhot)\T\mut$
is bounded in $\L\infty H${\ppg ,}
due to \eqref{primastima} and \eqref{quartastima}.
Therefore, by a standard argument
(multiply formally by $\xit$), we deduce that 
both $-\Delta\rhot$ and $\xit$ are bounded in the same space,
whence by elliptic regularity 
\Beq
  \norma\rhot_{\L\infty W} \leq c
  \aand
  \norma\xit_{\L\infty H} \leq c{\ppg ;}
  \label{daquarta}
\Eeq
{\ppg moreover,}
\Beq
  \norma\rhot_{\LQ\infty}
  + \norma{\psi(\rhot))}_{\LQ\infty}
  \leq c
  \quad \hbox{with $\psi=g,g',\pi$}{\ppg ,}
  \label{daquartabis}
\Eeq
due to the continuous embedding $W\subset\Linfty$
and the continuity of such~$\psi$'s.

\step
Fifth a priori estimate

{\ppg To prove an $L^\infty$ estimate 
rather than just a boundedness property, we borrow} the argument {\ppg in} \cite{CGPS3}.
We observe that the approximating solution {\ppg satisfies:} 
\Beq
  \frac 12 \, \dt \bigl[ \bigl( \coeff \bigr) \, |\mumkp|^2 \bigr]
  = \bigl[ \bigl( \coeff \bigr) \dt\mu + (\mu-k) g'(\rho) \dt\rho \bigr] \mumkp \,,
\Eeq
{\ppg for every $k\in\erre$}. 
Hence, by assuming {\ppg that }$k\geq\muz^*:=\norma\uz_\infty$
and {\ppg by} testing \eqref{primatau} {\ppg with} $\mumkp$, we obtain{\ppg :}
\Beq
  \frac 12 \iO (\coefft) |(\mu(t)-k)^+|^2
  + \intQt \kamurho |\nabla(\mu-k)^+|^2
  = - k \intQt \dt g(\rho) \, (\mu-k)^+{\ppg ;}
  \non
\Eeq
{\ppg with this, on} recalling that $g\geq0$ and $\kappa\geq\kmin$, we deduce the inequality
\Beq
  \frac 12 \iO |(\mu(t)-k)^+|^2
  + \kmin \intQt |\nabla(\mu-k)^+|^2
  = - k \intQt \dt g(\rho) \, (\mu-k)^+.
  \non
\Eeq
In~\cite{CGPS3}{\ppg , for} $\eps=1$, 
we have {\ppg that} $g(r)=r$ and $\kappa=1${\ppg ;}
the corresponding inequality is{\ppg :}
\Beq
  \frac 12 \iO |(\mu(t)-k)^+|^2
  + \intQt |\nabla(\mu-k)^+|^2
  = - k \intQt \dt\rho \, (\mu-k)^+.
  \label{daCGPStre}
\Eeq
Therefore, the argument used in {\ppg that} paper
can be repeated {\ppg here essentially} without changes{\ppg. As a matter of fact,}
the analogue of~\eqref{seconda} is never used {\ppg in \cite{CGPS3},}
the whole proof {\ppg being} based just on \eqref{daCGPStre},
the regularity $\dt\rho\in\L\infty H\cap\L2V$,
and an upper bound, say~$C_0$, for the corresponding norm{\ppg; moreover, the upper bound for $\mu$, that is constructed explicitely,
depends only} on {\pier $\Omega$, $T$,} $\muz^*$, and~$C_0$.
In the present case, we have to use the same regularity for~$\dt g(\rho)$
and estimate~\eqref{daquartabis}.
In conclusion, we~obtain{\ppg :}
\Beq
  \norma\mut_{\LQ\infty} \leq c.
  \label{quintastima}
\Eeq

\step
Sixth a priori estimate

We proceed formally, as done for the third a priori estimate,
by writing $\xi=\beta(\rho)$ as {\ppg if} $\beta$ were a smooth function
(see Remark~\ref{Noformal}).
We test {\ppg by~$(\xi(t))^5$} \eqref{secondatau}{\ppg ,} written in the form \eqref{perterza}{\ppg ,}
at (almost) any fixed time~$t\in(0,T)$.
We obtain{\ppg :}
\Beq
  5 \iO (\xi(t))^4 \beta'(\rho(t)) |\nabla\rho(t)|^2
  + \iO |\xi(t)|^6
  = \iO h(t) \, (\xi(t))^5 .
  \non
\Eeq
As the first term on the \lhs\ is nonnegative,
by the \holder\ inequality we deduce~that
\Beq
  \norma{\xi(t)}_6^6
  \leq \norma{h(t)}_6 \, \norma{(\xi(t))^5}_{6/5}
  = \norma{h(t)}_6 \, \norma{\xi(t)}_6^5{\ppg ,}
  \non
\Eeq
whence {\ppg he have} immediately {\ppg that} $\norma{\xi(t)}_6\leq\norma{h(t)}_6$.
We infer that
\Beq
  \norma{\Delta\rho(t)}_{\Lx6} \leq c \norma{h(t)}_{\Lx6}
  \aand
  \norma{\rho(t)}_{\Wx{2,6}} \leq c \norma{h(t)}_{\Lx6}{\ppg ,}
  \non
\Eeq
first by comparison in \eqref{perterza} 
and then by the standard regularity theory of linear elliptic equations.
As $\Wx{1,6}$ is continuously embedded in~$\Cx0$
(see~\eqref{sobolevbis}){\ppg , and as}
the above inequalities hold \aat,
we deduce that
\Beq
  \norma{\nabla\rho}_{\L2\Linfty} \leq c \norma h_{\L2{\Lx6}}.
  \non
\Eeq
Now, we observe that $h$ is bounded in $\L2{\Lx6}${\ppg ,}
thanks to \eqref{quartastima}, \eqref{daquartabis}, \eqref{quintastima},
and the Sobolev inequality.
Therefore, we conclude~that 
\Beq
  \norma{\nabla\rhot}_{\L2\Linfty} \leq c.
  \label{sestastima}
\Eeq
{\ppg A byproduct of our proof is that} $\norma\rhot_{{\pier \L2{\Wx{2,6}}}}\leq c$.

\step 
Seventh a priori estimate

{\ppg On recalling \accorpa{defbp}{defK}, the following preparatory identities hold for the approximating solution:} 
\Bsist
  && \nabla\Kmurho 
  = \kamurho \nabla\mu + \Kumurho \nabla\rho{\ppg ,}
  \label{auxnablaK}
  \\
  && \dt\Kmurho 
  = \kamurho \dt\mu + \Kumurho \dt\rho{\ppg ,}
  \label{auxdtK}
  \\
  && \dt\Kumurho 
  = \kapmurho \dt\mu + \Kdmurho \dt\rho .
  \label{auxdtKu}
\Esist
Moreover, we notice that 
\Beq
  \norma{\psi(\mu,\rho)}_{\LQ\infty}
  \leq c 
  \quad \hbox{with $\psi=\kappa,\kappa',K,\Ku, {\ppg \textrm{or~}}\Kd$}{\ppg ,}
  \label{daquinta}
\Eeq
{\ppg due to} our structural assumptions \eqref{hpkbis}
and {\ppg to} \eqref{quintastima}, \eqref{lingr}.
Now, we formally test \eqref{primatau} by $\dt\Kmurho$
and get{\ppg :}
\Bsist
  && \intQt (\coeff) \, \dt\mu \, \dt\Kmurho
  + \intQt \mu \, \dt g(\rho) \, \dt\Kmurho
  \non
  \\
  && {} + \intQt \kamurho \nabla\mu \cdot \nabla\dt\Kmurho
  = 0 .
  \label{testpersettima}
\Esist
{\pier With the help of} \accorpa{auxnablaK}{auxdtK}, {\pier we rewrite} the first two integral 
as follows{\ppg :}
\Bsist
  && \intQt (\coeff) \, \dt\mu \, \dt\Kmurho
  \non
  \\
  && = \intQt (\coeff) \, \kamurho \, |\dt\mu|^2
  + \intQt (\coeff) \, \Kumurho \, \dt\mu \, \dt\rho
  \non
  \\
  && \intQt \mu \, \dt g(\rho) \, \dt\Kmurho
  = \intQt \mu \, \dt g(\rho) \, \kamurho \, \dt\mu 
  + \intQt \mu \, \dt g(\rho) \, \Kumurho \, \dt\rho .
  \non
\Esist
In the third integral of \eqref{testpersettima}, 
we also integrate by parts and use \eqref{auxdtKu}.
We~get{\ppg :}
\Bsist
  && \intQt \kamurho \nabla\mu \cdot \nabla\dt\Kmurho
  = \intQt \bigl( \nabla\Kmurho - \Kumurho \nabla\rho \bigr) \cdot \nabla\dt\Kmurho
  \non
  \\
  && = \frac 12 \iO |\nabla\Kmurhot|^2
  - \iO \Kumurhot \nabla\rho(t) \cdot \nabla\Kmurhot
  - c
  \non
  \\
  && \quad {}
  + \intQt \nabla\Kmurho \cdot \dt \bigl( \Kumurho \nabla\rho \bigr) 
  \non
  \\
  && = \frac 12 \iO |\nabla\Kmurhot|^2
  - \iO \Kumurhot \nabla\rho(t) \cdot \nabla\Kmurhot
  - c
  \non
  \\
  && \quad {}
  + \intQt \nabla\Kmurho \cdot \nabla\rho \, \bigl( \kapmurho \dt\mu + \Kdmurho \dt\rho \bigr)
  + \intQt \Kumurho \nabla\Kmurho \cdot \nabla\dt\rho .
  \non
\Esist
{\ppg If we insert these} identities in \eqref{testpersettima},  {\ppg on}
keeping just the positive terms on the \lhs, {\ppg on}
recalling that $g\geq0$ and $\kappa\geq\kmin$, {\ppg and on}
and accounting for estimates \eqref{quintastima} and~\eqref{daquinta}, {\ppg then}
we deduce that
\Bsist
  && \kmin \intQt |\dt\mu|^2
  + \frac 12 \iO |\nabla\Kmurhot|^2
  \non
  \\
  && \leq c \intQt |\dt\mu| \, |\dt\rho|
  + c \intQt |\dt\rho|^2
  + c \iO |\nabla\rho(t)| \, |\nabla\Kmurhot|
  \non
  \\
  && \quad {}
  + c \intQt |\nabla\Kmurho| \, |\nabla\rho| \, \bigl( |\dt\mu| + |\dt\rho| \bigr)
  + c \intQt |\nabla\Kmurho| \, |\nabla\dt\rho|
  + c .
  \label{persettima}
\Esist
As the first three terms on the \rhs\ can be trivially dealt with 
by accounting for \eqref{primastima}, \eqref{secondastima}, 
and the elementary Young inequality, we {\ppg concentrate on} the last two integrals.
{\pier For every $\eps>0$, we deduce that}
\Bsist
  && \intQt |\nabla\Kmurho| \, |\nabla\rho| \, \bigl( |\dt\mu| + |\dt\rho| \bigr)
  + \intQt |\nabla\Kmurho| \, |\nabla\dt\rho|
  \non
  \\
  && \leq \eps \intQt |\dt\mu|^2
  + c_\eps \intQt |\nabla\rho|^2 \, |\nabla\Kmurho|^2
  \non
  \\
  && \quad {}
  + \intQt |\nabla\Kmurho|^2
  + \frac 12 \intQt |\nabla\rho|^2 \,|\dt\rho|^2 
  + \frac 12 \intQt |\nabla\dt\rho|^2 .
  \non
\Esist
On the other hand, we have{\ppg :}
\Bsist
  && \intQt |\nabla\rho|^2 \, |\nabla\Kmurho|^2 
  + \intQt |\nabla\Kmurho|^2
  + \intQt |\nabla\rho|^2 \,|\dt\rho|^2 
  + \intQt |\nabla\dt\rho|^2
  \non
  \\
  && \leq \iot \phi(s) \, \norma{\nabla\Kmurhos}_2^2 \ds
  + \iot \norma{\nabla\rho(s)}_\infty^2 \,\norma{\dt\rho(s)}_2^2  \ds
  + \intQt |\nabla\dt\rho|^2
  \non
  \\
  && \leq \iot \phi(s) \, \norma{\nabla\Kmurhos}_2^2 \ds
  + \norma{\nabla\rho}_{\L2\Linfty}^2 \norma{\dt\rho}_{\L\infty H}^2
  + \norma{\nabla\dt\rho}_{\L2H}^2\,,
  \non
\Esist
where $\phi(s):=\norma{\nabla\rho(s)}_\infty^2+1$.
As $\phi\in L^1(0,T)$ thanks to \eqref{sestastima}{\ppg , and as}
the last norms in the above inequality 
are bounded by \eqref{quartastima} and~\eqref{sestastima},
we can choose $\eps$ small enough and apply the Gronwall lemma.
We conclude~that
\Beq
  \norma{\dt\mut}_{\L2H}
  + \norma{K(\mut,\rhot)}_{\L\infty V} \leq c.
  \label{settimastima}
\Eeq

\step Consequence

By combining \eqref{settimastima}, \eqref{auxnablaK}, and $\kappa\geq\kmin$,
we derive~that
\Beq
  \norma{\nabla\mut}_{\L\infty H} \leq c,
  \quad \hbox{whence} \quad
  \norma\mut_{\L\infty V} \leq c.
  \label{stimamuLdueV}
\Eeq
Furthermore, by comparison in {\pier \eqref{primatau}}, we also deduce~that
\Beq
  \norma{\div(\kappa(\mut,\rhot)\nabla\mut)}_{\L2H} \leq c .
  \label{dasettima}
\Eeq

\step
Limit and conclusion

By the above estimates, there {\ppg are} a triplet $(\mu,\rho,\xi)$,
with $\mu\geq0$ \aeQ, and {\ppg two} functions $k$ and $\zeta$ such~that
\Bsist
  \hskip-1cm& \mut \to \mu
  & \quad \hbox{weakly star in $\H1H\cap\L2V\cap\LQ\infty$}{\ppg ,}
  \qquad
  \label{convmu}
  \\
  \hskip-1cm& \rhot \to \rho
  & \quad \hbox{weakly star in $\L\infty W$}{\ppg ,}
  \label{convrho}
  \\
  \hskip-1cm& \dt\rhot \to \dt\rho
  & \quad \hbox{weakly star in $\L\infty H\cap\L2V$}{\ppg ,}
  \label{convdtrho}
  \\
  \hskip-1cm& \xit \to \xi
  & \quad \hbox{weakly star in $\L\infty H$}{\ppg ,}
  \label{convxi}
  \\
  \hskip-1cm& K(\mut,\rhot) \to k  
  & \quad \hbox{weakly star in $\L\infty V$}{\ppg ,}
  \label{convKmu}
  \\
  \hskip-1cm& \div(\kappa(\mut,\rhot)\nabla\mut) \to \zeta
  & \quad \hbox{weakly in $\L2H$}{\ppg ,}
  \label{convdiv}
\Esist
at least for a susequence $\tau=\tau_i{\scriptstyle\searrow}0$.
By the weak convergence of time derivatives,
the Cauchy conditions \eqref{cauchy} hold for the limit pair $(\mu,\rho)$.
By \accorpa{convmu}{convdtrho}
and the compact embeddings $W\subset\Cx0$ and $V\subset H$,
we can apply well-known strong compactness results
(see, e.g., \cite[Sect.~8, Cor.~4]{Simon})
and{\pier , possibly taking another subsequence,} we have that
\Bsist
  & \mut \to \mu
  & \quad \hbox{strongly in $\L2H$ and \aeQ}
  \label{strongmu}
  \\
  & \rhot \to \rho
  & \quad \hbox{strongly in $\CQ$}
  \label{strongrho}.
\Esist
The weak convergence \eqref{convxi}{\pier , together with \eqref{strongrho},
implies} that $\xi\in\beta(\rho)$ \aeQ\ 
{\pier (see, e.g., \cite[Prop.~2.5, p.~27] {Brezis}), 
due to the maximal monotonicity of 
{\gianni the operator induced by $\beta$ on $L^2(Q)$}}.
Furthermore, the {\pier convergence stated in \accorpa{strongmu}{strongrho}
entails that} $\phi(\rhot)\to\phi(\rho)$ uniformly in $Q$
for $\phi=g,g',\pi$
and $\psi(\mut,\rhot)\to\psi(\mu,\rho)$ \aeQ\
for $\psi=\kappa,K$,
whence, in particular, $k=\Kmurho$.
As all the above functions $\psi(\mut,\rhot)$ 
are uniformly bounded,
we deduce that the convergence is in fact strong in $\LQ p$ for every $p<+\infty$
and weak star in $\LQ\infty$ in each case.
{\gianni This shows that the limits of the products 
\Beq
  (\coefftau) \, \dt\mut , \quad
  \mut \, g'(\rhot) \, \dt\rhot , \quad
  \kappa(\mut,\rhot) \, \nabla\mut ,
  \aand
  \mut \, g'(\rhot)
  \non
\Eeq
that appear in equations
\accorpa{primatau}{secondatau}
can be identified as the products of the corresponding limits. 
In particular, by using also \eqref{convdiv},
we derive that $\div(\kamurho\nabla\mu)$ equals $\zeta$ and belongs to~$\LQ2$.
All this implies both \eqref{prima} for the limit $(\mu,\rho)$
and the convergence of the normal trace $\kappa(\mut,\rhot)\,\nabla\mut\cdot\nu$.
Thus, the expected Neumann condition also holds in the limit{\ppg ,} and the proof is complete.}



\vspace{3truemm}

\Begin{thebibliography}{10}

\bibitem{Barbu}
V. Barbu,
``Nonlinear semigroups and differential equations in Banach spaces'',
Noord\-hoff,
Leyden,
1976.


\bibitem{Brezis}
H. Brezis,
``Op\'erateurs maximaux monotones et semi-groupes de contractions
dans les espaces de Hilbert'',
North-Holland Math. Stud.
{\bf 5},
North-Holland,
Amsterdam,
1973.


\bibitem{CGPS3} 
P. Colli, G. Gilardi, P. Podio-Guidugli, J. Sprekels,
Well-posedness and long-time behaviour for a nonstandard viscous Cahn-Hilliard 
system, {\giannipier {\it SIAM J. Appl. Math.} {\bf 71} (2011), 1849-1870.}

\bibitem{CGPSgen} 
P. Colli, G. Gilardi, P. Podio-Guidugli, J. Sprekels,
{\piecol Global existence and a new uniqueness proof for a singular/degenerate  
Cahn-Hilliard system with viscosity}, {\it manuscript.}





\bibitem{Podio}
P. Podio-Guidugli, 
Models of phase segregation and diffusion of atomic species on a lattice,
{\it Ric. Mat.} {\bf 55} (2006) 105-118.


\bibitem{Simon}
J. Simon,
{Compact sets in the space $L^p(0,T; B)$},
{\it Ann. Mat. Pura Appl.~(4)} {\bf 146} (1987) 65--96.

\End{thebibliography}

\End{document}

\bye